\documentclass[11pt]{amsart}
\usepackage[arrow,matrix,curve,cmtip,ps]{xy}
\newtheorem{thm}{Theorem}[section]   % Numbered within each section
\newtheorem{cor}[thm]{Corollary}     % Numbered along with thm
\newtheorem{lem}[thm]{Lemma}         % Numbered along with thm
\newtheorem{prop}[thm]{Proposition}  % Numbered along with thm

\theoremstyle{definition}
\newtheorem{defn}[thm]{Definition}   % Numbered along with thm

\newtheorem{rem}[thm]{Remark}        % Numbered along with thm
\newtheorem{ex}[thm]{Example}        % Numbered along with thm

\begin{document}

\title{$C^*$-equivalences of graphs}

\author{D. Drinen}
\email{Drinen@asu.edu}

\author{N. Sieben}
\email{nandor.sieben@asu.edu}
\address{Dept. of Mathematics\\Arizona State Univ.\\Tempe, AZ
85287}


\begin{abstract} 
Several relations on graphs, including primitive equivalence, explosion
equivalence and strong shift equivalence, are examined and
shown to preserve either the graph groupoid, a construction of Kumjian,
Pask, Raeburn, and Renault, or the groupoid of a pointed version of the
graph.  Thus these relations preserve either the isomorphism class or the
Morita equivalence class of the graph $C^*$-algebra, as defined by
Kumjian, Pask, and Raeburn. 
\end{abstract}

\maketitle

%------------------------------------------------------------------------
\section{Introduction}
\label{sec-intro}
%-----------------------------------------------------------------------

Given any finite square matrix $B$ with nonnegative integer entries and no
zero rows or columns, Cuntz and Krieger defined a $C^*$-algebra
$\mathcal{O}_B$, which is generated by partial isometries satisfying
relations associated to $B$ \cite{CK1}. Also, given any square matrix $B$
with nonnegative integer entries, we can build a directed graph by putting
$B_{ij}$ edges from vertex $i$ to vertex $j$. In \cite{kpr} Kumjian, Pask,
and Raeburn defined the graph $C^*$-algebra $C^*(E)$ of any countable
row-finite directed graph $E$ as the universal $C^*$-algebra generated by
a family of projections and partial isometries which satisfy relations
coming from $E$.  Given a graph $E$ with finitely many vertices and no
sources or sinks, the vertex matrix $B_E$ associated to the graph is
finite and has no zero rows or columns.  In this case it turns out that
$C^*(E)$ coincides with $\mathcal{O}_{B_E}$.  The graph thus becomes a
useful tool for visualizing and generalizing the Cuntz-Krieger algebras.

In \cite{kprr} Kumjian, Pask, Raeburn, and Renault defined the graph
groupoid $\mathcal{G}_E$ of any countable row-finite directed graph $E$
with no sinks.  In this case, the $C^*$-algebra of the groupoid coincides
with the $C^*$-algebra of the graph, so the graph groupoid is another tool
for understanding the Cuntz-Krieger algebras and a large class of the
graph algebras.  In this paper, we examine several equivalence relations
on graphs which preserve the graph groupoid or, in some cases, the
groupoid of a pointed version of the graph. 

In \cite{wat-prim} Enomoto, Fujii, and Watatani defined primitive
equivalence of finite directed graphs with no sources, no sinks, and no
multiple edges, and showed that it is a sufficient condition for
isomorphism of the resulting graph algebras.  In Section~\ref{sec-primeq}
we generalize their result to countable row-finite graphs.  Further,
Enomoto, Fujii, and Watatani showed that primitive equivalence is also a
necessary condition for isomorphism of the graph algebras of strongly
connected graphs with three vertices.  In
Section~\ref{sec-classification}, we show by counterexample that this does
not hold in the four-vertex case. 

Primitive equivalence involves changing the rows of a matrix.  In
graph-theoretical terms, this corresponds to changing the outgoing edges
at a vertex.  In Section~\ref{sec-revprimeq}, we define an equivalence
relation which involves changing the columns of the matrix (alternatively,
the incoming edges at a vertex).  We call this reverse primitive
equivalence, and show that it preserves the Morita equivalence class,
though not the isomorphism class, of the graph algebras.

Primitive equivalence and reverse primitive equivalence only make sense
for graphs with the same number of vertices.  In
Section~\ref{sec-explosions}, we define explosion and reverse explosion,
operations which can change the size of the graph.  The graph operation we
call reverse explosion was defined and called explosion in \cite{wat-prim}. 
We show that exploding a graph does not change the graph groupoid, hence
does not change the isomorphism class of its graph $C^*$-algebra.  Reverse
exploding a graph does not preserve the graph groupoid, but the resulting
graphs can be pointed in such a way that their groupoids are isomorphic. 
Thus, reverse exploding a graph preserves the Morita equivalence class of
its $C^*$-algebra.

In Section~\ref{sec-SSE}, we recall from \cite{ashton} the notion of
elementary strong shift equivalence of graphs, and show that elementary
strong shift equivalent graphs can be pointed in such a way that their
groupoids are isomorphic.  This is an alternate proof of a result of
Cuntz and Krieger in \cite{CK1}, which states that elementary strong
shift
equivalent matrices correspond to Morita equivalent Cuntz-Krieger 
algebras. We then examine
the relationship between elementary strong shift equivalence and explosion
equivalence. 

The authors would like to thank Alex Kumjian, David Pask, John Quigg, and
Jack Spielberg for many helpful discussions.

%=================================================================
\section{Preliminaries}
\label{sec-prelim}
%==================================================================

In \cite{wat-prim}, a $C^*$-algebra was associated to every connected 
finite directed
graph with no multiple edges, no sources, and no sinks.  We review this
construction as we set up the notation.  A {\it directed graph} $E$
consists of a set $E^0$ of vertices, a set $E^1$ of edges and maps
$s,r:E^1 \rightarrow E^0$ describing the source and range of each
edge.  Denote by $E^j$ the set of paths in $E$ of length $j$.  Here,
zero-length paths (i.e., vertices) are allowed.  Denote by $E^*$
the set of all finite paths in $E$ and by 
$E^\infty$ the infinite one-sided path space of 
$E$.  We extend $s$ and $r$ to $E^*$ and $s$ also to $E^\infty$. 
Associated to every directed graph $E$ is an $E^0 \times
E^0$ {\it vertex matrix} $B_E$, defined by $B_E(v,w) = \#\{e \in
E^1\,|\,s(e)
= v, r(e) = w\}$.  That is, the $(v,w)$ entry of $B_E$ is the number of
edges
in $E$ from $v$ to $w$. $E$ has no multiple edges if and
only
if $B_E$ is a 0-1 matrix.  $E$ has no sources (resp. sinks) if and only if
$B_E$ has no zero columns (resp. rows).  A directed graph is said to be
{\it strongly connected}
if for every pair of vertices $v$ and $w$, there is a path 
from $v$ to $w$ and a path from $w$ to $v$.  A directed graph is said to
be
{\it connected} if between each pair $v, w$
there is an undirected path from $v$ to $w$.  

In this paper, all graphs
are assumed to be countable, directed, connected, and to have no multiple
edges.

Note that in \cite{wat-prim}, Enomoto, Fujii and Watatani worked with
the adjacency matrix instead of the vertex matrix, whose transpose is the
adjacency matrix.  We choose to work with the vertex matrix in order
to
be consistent with the more recent graph algebra literature \cite{kprr}.

%Also note that $B_E$ satisfies condition (I) defined in \cite{CK1} if and
%only if every vertex of $E$ connects to a vertex which has two distinct
%elementary return paths.  In \cite{kpr}, it is shown that, if $E$ is
%finite and has no sinks, condition (I) for $B_E$ is equivalent to the
%condition that every loop in $E$ has an exit, which is referred to as
%condition (L).  A finite 0-1 matrix $B$ is said to be {\it irreducible}
%if for
%every $i,j$, there exists an $N$ such that $B^N(i,j) > 0$.  Note that a 
%graph is strongly connected if and only if its vertex matrix is
%irreducible, and that any
%irreducible matrix which is not a permutation matrix satisfies (I).

For any row-finite graph $E$, a Cuntz-Krieger $E$-family is a set 
$\{P_v\,|\,v \in E^0\}$ of mutually orthogonal projections together with
a set $\{S_e\,|\,e \in E^1\}$ of partial isometries satisfying:
\begin{enumerate}
\item[(a)] $S_e^*S_e = P_{r(e)}$;
\item[(b)] $P_v = \sum_{s(e)=v}S_eS_e^*$, for $v \in s(E^1)$.
\end{enumerate}
Kumjian, Pask, and Raeburn \cite{kpr} defined the $C^*$-algebra of the
graph, denoted by $C^*(E)$, to be the universal $C^*$-algebra generated by
a Cuntz-Krieger $E$-family.

We now recall the construction of a groupoid from a
row-finite graph \cite{kprr}.  For $x,y \in E^{\infty},\,k \in
{\bf Z},$ say
$x \sim_k y$ if and only if $x_i = y_{i-k}$ for large $i$, where $x_i$
denotes the $i$-th edge of $x$.  We remark that this definition differs 
slightly from the one given in \cite{kprr}, but it coincides with the 
currently accepted convention.
Define $\mathcal{G}_E,$ the path groupoid
of $E,$ by $\mathcal{G}_E = \{(x,k,y)\,|\,x \sim_k y\}.$  The groupoid
operations in $\mathcal{G}_E$ are
$$
(x,k,y)^{-1} = (y,-k,x) \qquad \hbox{ and } \qquad (x,k,y) \cdot (y,l,z)
= (x,k+l,z).
$$

Alternatively, one can define $$\mathcal{G}_E = \{(\alpha, x, \beta) \in
E^* \times E^\infty \times E^*\,|\,r(\alpha) = r(\beta) = s(x)\}/\sim,$$

\noindent where $\sim$ denotes the equivalence relation $(\alpha, \gamma
x, \beta) \sim (\alpha \gamma , x, \beta \gamma)$.  To see that the two
definitions coincide, the reader may check that the map $[\alpha, x,
\beta] \mapsto (\alpha x, |\alpha| - |\beta|, \beta x)$ is a groupoid
isomorphism.  With
this definition, we find that $[\alpha, x, \beta]^{-1} = [\beta, x,
\alpha]$ , that $[\alpha, x, \beta]$ and $[\gamma, y,
\delta]$ are composable if and only if $\beta x = \gamma y$, and that 
\[ [\alpha, x, \beta] [\gamma, y, \delta] = \left\{
\begin{array}{ll}
   [\alpha, x, \delta \epsilon] & \hbox{if $\gamma \epsilon = \beta$ and
$y = \epsilon x$} \\
 
[\alpha \epsilon , y, \delta] & \hbox{if $\beta
\epsilon =
\gamma$ and $x = \epsilon y$}
\end{array}
\right.  \hbox{.} \]

With the topology generated by the sets $Z(\alpha,
\beta) := \{[\alpha, x, \beta]\,|\,s(x) = r(\alpha)\}$, $\mathcal{G}_E$
is a locally compact Hausdorff $r$-discrete groupoid with Haar system.
If $E$ has
no sinks, then $C^*(E)$, the $C^*$-algebra of $E$ constructed in
\cite{kpr}, coincides with $C^*(\mathcal{G}_E)$.  

We now seek to remove the restriction on sinks.  First recall from \cite{kprr}
that a pair $(E,S)$, where $E$ is a row-finite graph  with no sinks and $S$ is
a set of vertices of $E$, is called a {\it pointed graph}.   If $(E,S)$ is a
pointed graph, then $S$ determines a clopen subset $\{x \in E^\infty\,|\,s(x)
\in S\}$ of the unit space of $\mathcal{G}_E$, which we also denote by $S$.  If
$S$ is {\it cofinal}, meaning that given any $x \in E^\infty$ there exists $v
\in S$ and a finite path from $v$ to $s(x_i)$ for some $i$, then
$C^*(\mathcal{G}_{(E,S)})$ is Morita equivalent to
$C^*(\mathcal{G}_E)$, where $\mathcal{G}_{(E,S)}$ denotes
$\mathcal{G}_E$ restricted to $S$ \cite{kprr}.

Given a row-finite graph $E$ with sinks at
$\{v_i\}_{i \in I}$,  
define a pointed graph $\tilde{E}$ as follows:  the vertices of
$\tilde{E}$ are
the
vertices of $E$, along with the additional vertices $\{w^j_i\}$ for $i \in
I$, $j=1,2,\ldots$.  The edges in $\tilde{E}$ are the edges in $E$ along
with
an
edge from $v_i$ to $w_i^1$ for every $i \in I$, and an edge from $w_i^j$
to $w_i^{j+1}$ for every $i \in I$, $j = 1,2,\ldots$.  We have simply
added a distinct infinite tail to each sink in $E$.  Define the pointing
set of $\tilde{E}$ to be the original set of vertices of $E$.  The reader
may verify that $C^*(\mathcal{G}_{\tilde{E}}) \cong C^*(E)$ by checking
that both are generated by the same family of projections and partial
isometries.  

%------------------------------------------------------------------------
\section{Primitive Equivalence}
\label{sec-primeq}
%------------------------------------------------------------------------

The following definition is due to Enomoto, Fujii, and Watatani.  Let $B$
be an $n \times n$, 0-1 matrix. For $1 \leq p \leq n$, denote
the $p$-th row of $B$ by $B_p$ and denote the row
$(0,\dots,0,1,0,\dots,0)$, where the $1$ is in the $i^{th}$ position, by
$E_i$.  
We apologize for any confusion caused by this
multiple use of the letter $E$, but we are using standard conventions of
\cite{wat-prim} for primitive transfer and standard conventions of \cite{kpr} 
for graphs.

Now suppose that there is a $p$ such that $B_p$ is not a zero row
and  
$$B_p = E_{k_1} + \cdots + E_{k_r} + B_{m_1} + \cdots + B_{m_s}$$ 

\noindent for some distinct $k_1, \dots ,k_r$, $m_1, \dots ,
m_s$ such that $p \not \in \{m_1, \dots, m_s\}$ and $B_{m_i}$ is not a zero
row for any $i$.  
Define a new matrix $C$
by 
$$
C_{ij} = \left\{  \begin{array}{ll}
			B_{ij} & \hbox{if $i \neq p$} \\
			1      & \hbox{if $i = p$ and $j \in \{k_1,
\dots, k_r, m_1, \dots, m_s\}$} \\
			0      & \hbox{otherwise.}
		    \end{array}
	   \right.  $$

\noindent That is, start with B, zero out the $p$-th row, and then put
back $1$'s
in positions $k_1, \dots, k_r, m_1, \dots, m_s$.  $C$ is called a {\it
primitive transfer} of $B$ at $p$.  Note that this notion does not depend
on finiteness of the matrix $B$.

\begin{defn}
If $B$ and $C$ are 0-1, square matrices of the same (possibly infinite) size, 
we say $B$ is {\it
primitively equivalent} to $C$ if and only if there exist matrices $D_1,
\dots, D_q$ such that $D_1 = B$, $D_q = C$, and for every $1 \leq i \leq
q-1$, one of the following holds:

\begin{enumerate}
\item[(i)] $D_i$ is a primitive transfer of $D_{i+1};$
\item[(ii)] $D_{i+1}$ is a primitive transfer of $D_i;$
\item[(iii)] $D_i = PD_{i+1}P^{-1}$ for some permutation matrix $P.$
\end{enumerate}
\end{defn}

We say that two matrices which satisfy the third condition above are
{\it permutations} of each other.  Matrices which are permutations of
each
other should be primitively equivalent because we would like the vertex 
matrices
of isormorphic graphs to be primitively equivalent.  This is implicit but never 
stated in \cite{wat-prim}.  There are matrices
which cannot be primitively transferred in any number of steps into a
permutation of the same matrix.  For example, 
$$A = \left( \begin{array}{ccc}
		1 & 1 & 1\\
		1 & 0 & 1\\
                1 & 0 & 0
	     \end{array}
      \right) \hbox{ and }
  B = \left( \begin{array}{ccc}
		0 & 1 & 1\\
		1 & 1 & 1\\
                0 & 1 & 0
	     \end{array}
      \right) \hbox{,}   
$$              
\noindent are permuations of each other, but the reader may check that
they would not be primitively equivalent if condition (iii) were removed
from the definition. 

Franks \cite[Corollary 2.2]{franks} defined a similar matrix operation.
His operation applies to multigraphs and it involves only two rows or
columns.  He used this move in finding a canonical form for the flow
equivalence class of a matrix.

The primitive transfer has the following graph-theoretical interpretation:
suppose 
vertex $p$ points to the same vertices which are pointed to by vertices $m_1,
\dots, m_s$ (and only one of the vertices $m_1, \dots , m_s$ points to each
of those vertices) and, in addition, vertex $p$ points to vertices $k_1,
\dots, k_r$.  Then a {\it primitively transferred graph} can be obtained
by
erasing all the edges emanating from vertex $p$, except those pointing to
vertices $k_1, \dots, k_r$, and adding an edge from vertex $p$ to each of the
vertices $m_1, \dots, m_s$. Note that this procedure is only allowed if we do
not create any multiple edges. The following example shows that we may first
erase and then redraw the same edge.  

\begin{ex} \label{ex-prim} Consider the following graph $E$ and its primitive
transfer $F$: 
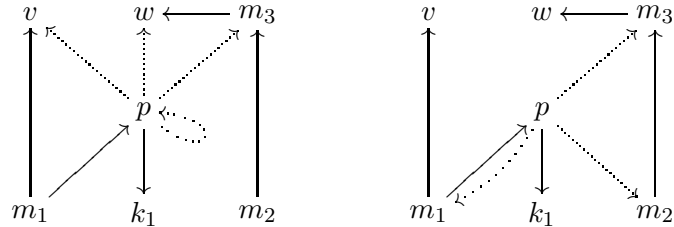
$$ 
\xymatrix{  
v	& w	& m_3 \ar[l]						\\
	& p \ar[d] \ar@{.>}[u] \ar@{.>}[lu]  \ar@{.>}[ru] \ar@{.>}@(dr,r)[] \\ 
m_1 \ar[uu] \ar[ur]	& k_1	& m_2\ar[uu]  
} 
\qquad\qquad 
\xymatrix{ 
v	&w	& m_3 \ar[l]            			\\
	& p \ar[d] \ar@{.>}@/^/[ld] \ar@{.>}[rd] \ar@{.>}[ur]	\\ 
m_1 \ar[uu]\ar[ur]  & k_1                 & m_2 \ar[uu]  
} 
$$ 
\end{ex}

\noindent Vertex $p$ points to vertices $p$, $v$, $w$, $m_3$ and $k_1$.
Together, vertices $m_1$, $m_2$ and $m_3$ point to vertices $p$, $v$, $w$ and
$m_3$. This means that, if $B$ is the vertex matrix of $E$, then $B_p =
E_{k_1} + B_{m_1} + B_{m_2} + B_{m_3}$.  So row $p$  of the vertex matrix
of $F$ has 1's in positions $k_1$, $m_1$, $m_2$  and $m_3$.  Thus $F$ has
edges from vertex $p$ to vertices $k_1$, $m_1$, $m_2$ and $m_3$.  

\begin{defn}
Two graphs with no multiple edges are said to be {\it primitively
equivalent} if and only if their vertex matrices are primitively
equivalent. 
\end{defn}

We will show that if $F$ is a primitive transfer of $E$,
then $\mathcal{G}_E$
and $\mathcal{G}_F$ are isomorphic.  But first we need the following
definitions, notation and lemmas.

Let $E$ and $F$ be row-finite graphs with no sinks. To simplify notation, we
denote
their vertex matrices by $B$ and $C$, respectively.  Suppose that $F$ is a
primitive tranfer of $E$.  Without loss of generality, assume that $1 \in E^0$
and $B_1 = E_{k_1}
+ \cdots + E_{k_r} + B_{m_1} + \cdots + B_{m_s}$.  We identify $E^0$ and
$F^0$.  Define $K := \{k_1,\dots, k_r\}$ and $M := \{m_1, \dots, m_r\}$.
Note that $K \cap M = \emptyset$.  Since $E$ and
$F$ have no multiple edges we can use the notation $e^{ij}$ to denote the
unique edge in $E$ with source $v_i$ and range $v_j$, if there is one,
that is, if $B_{ij} = 1$.  Likewise, we will denote edges in $F$ by
$f^{ij}$.

The following lemma is an immediate consequence from the definition of the
primitive transfer.

\begin{lem}
\label{lem-KUM}
If $f^{1j}$ exists, then $j \in K \cup M.$
\end{lem}

If $f \in F^1$ is of the form $f^{1m}$ for some $m \in M$, then we will
call $f$ a {\it new} (or {\it newly introduced}) edge.  

\begin{lem}
\label{lem-fijnotneweijexists}
If $f^{ij}$ is not new, then $e^{ij}$ exists.
\end{lem}

\begin{proof}
$C_{ij} = 1$ because $f^{ij}$ exists.  Now, if $f^{ij}$ is not
new, then either $i \neq 1$ or $j \not\in M$. If $i \neq 1$, then the
primitive transfer does not change row $i$, so $B_{ij} = C_{ij} = 1$.  If,
on the other hand, $i = 1$ and $j \not\in M$, then $j \in K$ by
Lemma~\ref{lem-KUM}.  So $B_i = \cdots + E_j + \cdots$ and thus $B_{ij} =
1$.  In either case, $B_{ij} = 1$, so $e^{ij}$ exists.
\end{proof}

\begin{lem}
\label{lem-no2consecnew}
If $f$ and $f'$ are consecutive edges in $F$ {\rm (}that is, $r(f) =
s(f')${\rm )}, then
$f$ and $f'$ cannot both be new.
\end{lem}

\begin{proof}
This follows from the definition of new and the fact that $1 \not\in M$.
\end{proof}

\begin{lem}
\label{lem-newoldhasinverse}
If $f^{1m}f^{mj} \in F^2$ for some $m \in M$, then $e^{1j}$ exists, and $j
\not \in K$.
\end{lem}

\begin{proof}
Since $f^{mj}$ is not new, $e^{mj}$
exists by Lemma~\ref{lem-fijnotneweijexists}, and so $B_{mj} = 1$.  Hence
we have $B_{1j} = 1$, since $B_1 = \cdots + B_m + \cdots$.  Thus $e^{1j}$
exists.  

Now suppose $j \in K$.  Then we have $B_1 = \cdots + E_j + \cdots + B_m +
\cdots$.   But we know that $B_{mj} = 1$, and thus $B_{1j} > 1$, a
contradiction.
\end{proof}

\begin{prop}
\label{prop-main}
If $E$ is a row-finite graph, and $F$ is a primitive
transfer of $E$, then $\mathcal{G}_E \cong \mathcal{G}_F$.
\end{prop}

\begin{proof} The strategy of the proof is as follows:  We use the
properties
of the primitive transfer to construct an $s,r$-preserving injective map $\phi$
from the edges in $E$ to paths of length one or two in $F$.  This will
induce
injective  maps from finite paths in $E$ to finite paths in $F$, from infinite
paths in $E$ to infinite paths in $F$, and from $\mathcal{G}_E$ to
$\mathcal{G}_F$.  The injective map between the groupoids will turn out to be a
surjective homomorphism.

We use the same notation as above and assume, without loss of generality,
that $B_1 = E_{k_1} + \cdots + E_{k_r} + B_{m_1} + \cdots + B_{m_s}$. 

Note that if $E$ has an edge from vertex 1 to vertex $j$ for some $j
\not\in K$,
then there is a unique $m \in M$ such that
$B_{mj} = 1$.  Since $C_{1m} = 1$, we can  define $\phi: E^1 \rightarrow
F^1
\cup F^2$ by
$$
\phi(e^{ij}) = \left\{ \begin{array}{ll}
       f^{1m}f^{mj} & \hbox{if $i = 1$ and $j \not\in K$}
\\
       f^{ij}       & \hbox{else.}
		       \end{array}
	       \right.
$$

For instance, in Example \ref{ex-prim}, $\phi(e^{pp}) = f^{pm_1}f^{m_1p}$,
$\phi(e^{pv}) = f^{pm_1}f^{m_1v}$, $\phi(e^{pw}) = f^{pm_3}f^{m_3w}$ and
$\phi(e^{pm_3})=f^{pm_2}f^{m_2m_3}$. All the other edges would be mapped
to their corresponding edges.

$\phi$ induces a map, which we also denote by $\phi$, from $E^* \rightarrow
F^*$ by 
$$
\phi(\alpha_1\alpha_2\cdots\alpha_{|\alpha|}) =
 \phi(\alpha_1) \phi(\alpha_2) \cdots \phi(\alpha_{|\alpha|}).
$$
$\phi:E^\infty \rightarrow F^\infty$ is defined similarly, and
$\phi:\mathcal{G}_E \rightarrow \mathcal{G}_F$ is defined by 
$$
\phi [\alpha, x,\beta] = [\phi(\alpha), \phi(x), \phi(\beta)].
$$  
It is easily seen that $\phi:\mathcal{G}_E \rightarrow \mathcal{G}_F$ is
a well-defined homomorphism. 

In order to show that $\phi:\mathcal{G}_E \rightarrow \mathcal{G}_F$ is
injective, we need to know that $\phi:E^* \rightarrow F^*$ and $\phi:
E^\infty \rightarrow F^ \infty$ are injective.  

First note that $K$ and $M$ are disjoint. Hence if
$\phi(e)_1=\phi(e')_1$ for some $e,e'\in E^1$ then $|\phi(e)|=|\phi(e')|$.  
Now suppose that $\phi(\alpha) = \phi(\beta)$ for some finite or infinite paths
$\alpha$ and $\beta$. It follows by induction that
$\phi(\alpha_i)=\phi(\beta_i)$ for all $i$. Since $\phi: E^1 \rightarrow F^1
\cup F^2$ is clearly injective, we have that $\phi: E^* \rightarrow F^*$ and
$\phi: E^\infty \rightarrow F^\infty$ are injective.

We are now in position to show that
$\phi:\mathcal{G}_E \rightarrow \mathcal{G}_F$ is injective.  Suppose  
$\phi[\alpha, x, \beta] = \phi[\gamma, y, \delta]$.  Then we can assume,
without loss 
of generality, that $\phi(\alpha) = \phi(\gamma) \eta$, $\phi(\beta) =
\phi(\delta) \eta$, and $\phi(y) = \eta \phi(x)$.  We claim that $\eta \in
\phi(E^*)$.  If $|\phi(\alpha)| \leq 1$ then either $\eta = \phi(\alpha)$ or $\eta =
r(\phi(\alpha)) = \phi(r(\alpha))$, so we can suppose that   
$|\phi(\alpha)| > 1$.  Since $\phi(\gamma) \eta =
\phi(\alpha_1) \cdots \phi(\alpha_{|\alpha|})$ and $\phi(\alpha_i)$ has
length one or two for every $i$, it follows that for some $k$, either
$\phi(\gamma) =   
\phi(\alpha_1) \dots \phi(\alpha_k)$ or $\phi(\gamma) = 
\phi(\alpha_1) \dots \phi(\alpha_k) f$, where $\phi(\alpha_{k+1}) = fg$
($f,g \in F^1$).  However, the latter case is not possible since, by
definition of $\phi$, the last edge of $\phi(\gamma)$ cannot be a new
edge, but $f$ must be a new edge.  In the former case, $\eta =
\phi(\alpha_{k+1} \dots \alpha_{|\alpha|})$.  Thus there exists $\mu \in
E^*$ such that $\eta = \phi(\mu)$.  So $\phi(\alpha) = \phi(\gamma) \phi(\mu) = \phi(\gamma \mu)$.  Similarly,
$\phi(\beta) = \phi(\delta \mu)$, and $\phi(y) = \phi(\mu x)$.
By injectivity of $\phi:E^* \rightarrow F^*$ and $\phi: E^\infty
\rightarrow F^\infty$, we have $\alpha = \gamma \mu$, $\beta = \delta
\mu$, and $y = \mu x$, and hence $[\alpha, x, \beta] = [\gamma, y,
\delta]$.  Thus $\phi: \mathcal{G}_E \rightarrow \mathcal{G}_F$ is
injective.  

We now show that $\phi:\mathcal{G}_E \rightarrow \mathcal{G}_F$
is onto.  Since elements of the form $[f, y, r(f)]$, where $f \in F^1$, and $y
\in F^\infty$ with $r(f) = s(y)$, generate $\mathcal{G}_F$, it suffices to show
that each of them is in the range of $\phi$.

Use the following procedure to find an inverse image for $y$.  Note that
every new edge is followed by an edge which is not new
(Lemma~\ref{lem-no2consecnew}).  Lemma~\ref{lem-newoldhasinverse} says
that we can find an inverse image for these new-not new pairs.  The
remaining edges are
not new, and by Lemma~\ref{lem-fijnotneweijexists} these edges can be
pulled
back individually. 

Now we fix an edge $f = f^{ij} \in F^1$ and a path $y \in
F^\infty$ with $s(y) = r(f)$.  If $f$ is not new, then clearly
$\phi[e^{ij},
\phi^{-1}(y), r(e^{ij})] = [f, y, r(f)]$.
If, on the other hand, $f = f^{1m}$ is a new edge, then there
must be an edge $e' \in E^1$ such that $\phi(e') = fy_1$.
Further, if $f$ is a new edge, then $y_1$ cannot be new, so there must be
an edge $e'' \in E^1$ with $\phi(e'') = y_1$. In
this case, 
$$\phi[e', \phi^{-1}(y_2y_3\cdots), e''] = [fy_1, y_2y_3\cdots, y_1] = 
[f, y, r(f)].$$
Thus $\phi$ is onto. 

All that remains to show is that $\phi$ is continuous and open.  To see
that $\phi$ is open, the reader may check that $\phi(Z(\alpha, \beta)) = 
Z(\phi(\alpha), \phi(\beta))$ for any finite paths $\alpha$ and $\beta$.  
Likewise, $\phi^{-1}(Z(\gamma, \delta)) = Z(\phi^{-1}(\gamma),
\phi^{-1}(\delta))$ for any $\gamma, \delta$, so $\phi$ is continuous.
\end{proof}

The proposition immediately yields the following corollary, which was
proved
in \cite{wat-prim} for the case where $E$ and $F$ are finite graphs which 
satisfy (L) and have no sources or sinks:

\begin{cor}
\label{cor-maincor}
If $E$ is a row-finite graph with no sinks and $F$ is primitively
equivalent to $E$, then $\mathcal{G}_E \cong \mathcal{G}_F$
and hence $C^*(E) \cong C^*(F)$. 
\end{cor}

Recall that if $E$ has sinks, we can build a graph $\tilde{E}$ with no
sinks by
affixing a distinct infinite tail to each sink.  By pointing $\tilde{E}$
at all the original vertices of $E$, we obtain a pointed graph whose
groupoid $C^*$-algebra coincides with $C^*(E)$.

\begin{lem}
\label{lem-sinks}
Let $E$ be a row-finite graph, possibly with sinks, and let $F$ be a
primitive
transfer of $E$.  Then $\tilde{F}$ is a primitive transfer of $\tilde{E}$.    
\end{lem}

\begin{proof}
Let $B$, $C$, $\tilde{B}$, and $\tilde{C}$ denote the vertex matrices of
$E$, $F$, $\tilde{E}$, and $\tilde{F}$, respectively.  Assume that $B_p =
E_{k_1} + \cdots E_{k_r} + B_{m_1} + \cdots + B_{m_s}$, so that $C_p =
E_{k_1} + \cdots E_{k_r} + E_{m_1} + \cdots + E_{m_s}$.  Since for each
$i$, vertex $m_i$ is not a sink in $E$, $B_{m_i j} = 1$ if and only if
$\tilde{B}_{m_i j} = 1$.  Similarly, since vertex $p$ is not a sink in
$E$, $B_{pj} = 1$ if and only if $\tilde{B}_{pj} = 1$, Thus $\tilde{B}_{p}
= E_{k_1} + \cdots E_{k_r} + \tilde{B}_{m_1} + \cdots + \tilde{B}_{m_s}$,
and $\tilde{C}_p = E_{k_1} + \cdots E_{k_r} + E_{m_1} + \cdots + E_{m_s}$.
Thus $\tilde{F}$ is a primitive transfer of $\tilde{E}$.   
\end{proof}

\begin{thm}
If $F$ is primitively equivalent to the row-finite graph $E$,
then $C^*(E) \cong C^*(F)$.  
\end{thm}

\begin{proof}
By Lemma~\ref{lem-sinks} and
Corollary~\ref{cor-maincor} we have $C^*(\mathcal{G}_{\tilde{E}}) \cong
C^*(\mathcal{G}_{\tilde{F}})$, and hence $C^*(E) \cong C^*(F)$.   
\end{proof}

%----------------------------------------------------------------------
\section{Classification}
\label{sec-classification}
%----------------------------------------------------------------------

A matrix $B$ is said to be {\it irreducible} if for every $i, j$, there
exists an $N \in {\bf N}$ such that $B^N(i,j) > 0$.  In \cite{wat-prim} a
computer, along with some $K$-theory, was used to show that for all $3
\times 3$ irreducible matrices which are not permutation matrices,
primitive equivalence is a necessary as well as sufficient condition for
isomorphism of the Cuntz-Krieger algebras. 

We used a similar method to see whether this result is true for
irreducible $4 \times 4$ matrices. It is not.  The following
counterexample is one of many: 
$$A = \left( \begin{array}{cccc}
		1 & 1 & 1 & 1 \\
		1 & 0 & 1 & 1 \\
		1 & 1 & 0 & 1 \\
		1 & 1 & 1 & 0
	     \end{array}
      \right) \hbox{ and }
B = \left( \begin{array}{cccc}
		0 & 1 & 1 & 1 \\
		1 & 0 & 1 & 1 \\
		1 & 1 & 0 & 1 \\
		1 & 1 & 0 & 0
	     \end{array}
      \right) 
$$
By \cite{pask/rae-K}, $K_0(\mathcal{O}_A)$ is the abelian group generated
by $\{[P^A_i]\,|\,i = 1, \dots , 4\}$ subject to the relations $[P^A_i] =
\sum_j
A(i,j)[P^A_j]$, and similarly for $\mathcal{O}_B$.  One may check that 
$$[P^A_1] \mapsto (0,2), [P^A_2] \mapsto (0,1), [P^A_3] \mapsto (1,4),
[P^A_4] \mapsto (1,1)$$ 
\noindent is a faithful representation of
$K_0(\mathcal{O}_A)$ as ${\bf Z}_2 \oplus {\bf Z}_6$, and that 
$$
[P^B_1] \mapsto (0,1), [P^B_2] \mapsto (1,1), [P^B_3] \mapsto
(0,4), [P^B_4] \mapsto (1,2)$$
\noindent is a faithful representation of
$K_0(\mathcal{O}_B)$ as ${\bf Z}_2 \oplus {\bf Z}_6$.  Further, notice
that $[1_{\mathcal{O}_A}] = \sum_1^4 [P^A_i] \mapsto (0,2)$, and
$[1_{\mathcal{O}_B}] \mapsto (0,2)$, as well.  Thus, since $A$ and
$B$ are irreducible and there is an
isomorphism between the $K_0$ groups which preserves the class of the
identity, $\mathcal{O}_A \cong \mathcal{O}_B$ by \cite{rordam-K}. 

However, it is easy to check by hand using the definition of the primitive
transfer that each of these matrices is primitively equivalent only to its
permutations, and that $A$ and $B$ are not permutations of each other.

%======================================================================
%
\section{Reverse Primitive Equivalence}
\label{sec-revprimeq}
%======================================================================

In this section, we define a modified version of primitive
equivalence using column operations instead of row operations.  Recall
that a cofinal vertex is one from which any infinite path can be
intercepted.

\begin{defn}
Suppose that $B$ and $C$ are 0-1 (possibly infinite) square matrices.  We
say $C$ is a {\it reverse primitive transfer} of $B$ if $C^T$ is a
primitive transfer of $B^T$ at a cofinal vertex.  We say that $B$ and $C$
are {\it reverse primitively equivalent} if there is a sequence $B = D_1,
D_2, \dots, D_q = C$ such that for each $i < q$, $D_{i+1}$ is a
permutation of $D_i$, $D_{i+1}$ is a reverse primitive transfer of $D_i$,
or $D_i$ is a reverse primitive transfer of $D_{i+1}$.
\end{defn}

Two graphs $E$ and $F$ are \emph{reverse} graphs if the vertex matrices
$B_E$ and $B_F$ are transposes of each other, that is, $E$ and $F$ have
the same vertices and their edges have opposite directions.  Disregarding 
cofinality, two graphs are reverse primitive equivalent if their reverse
graphs are primitively equivalent. The following example shows that
reverse primitive equivalence of graphs $E$ and $F$ does not imply that
$C^*(E) \cong C^*(F)$. 

\begin{ex}
\label{ex-rpe}
If 
$$
B = \left( \begin{array}{cccc}
                0 & 1 & 1 & 1 \\
                1 & 0 & 0 & 0 \\
                1 & 0 & 0 & 0 \\
                1 & 0 & 0 & 0
             \end{array}
      \right) \hbox{ and }
C = \left( \begin{array}{cccc}
                0 & 1 & 0 & 1 \\
                1 & 0 & 1 & 0 \\
                1 & 0 & 0 & 0 \\
                1 & 0 & 0 & 0
             \end{array}
      \right) \hbox{,}
$$
\noindent then $C^T$ is a primitive transfer of $B^T$ (via $B^T_3 =
B^T_2$).  However, $\mathcal{O}_C \cong \mathcal{O}_3$ and  $\mathcal{O}_B
\cong \mathcal{O}_3 \otimes M_2$, so they are not isomorphic
\cite{paschke-salinas}. 
\end{ex}

Suppose that $E$ and $F$ are row-finite graphs and that $F$ is a
reverse primitive transfer of $E$ at vertex 1.  That is, $B_1^T = E_{k_1}
+ \cdots + E_{k_r} + B_{m_1}^T + \cdots + B_{m_s}^T$.  Where appropriate,
we use the notation established prior to Proposition~\ref{lem-KUM}.  

We have analogues of some of the preliminary lemmas for
Proposition~\ref{prop-main}, and their proofs are similar: 

\begin{lem}
\label{lem-KUM2}
If $f^{i1}$ exists, then $i \in K \cup M.$
\end{lem}

The definition of {\it new} must be altered slightly.  $f \in F^1$ is
said to be a {\it new} edge if it is of the form $f^{m1}$ for some $m \in
M$.

\begin{lem}
\label{lem-oldnewhasinverse}
If $f^{im}f^{m1} \in F^2$ for some $m \in M$, then $e^{i1}$ exists.
\end{lem}

Note that, with the modified definition of new,
Lemma~\ref{lem-fijnotneweijexists} and Lemma~\ref{lem-no2consecnew} are
true exactly as stated.

\begin{prop}
\label{prop-reverse}
If the row finite  graph $F$ is a reverse primitive
transfer of $E$ at $v$, then
$\mathcal{G}_{(E,\{v\})}
\cong
\mathcal{G}_{(F,\{v\})}$.
\end{prop}

\begin{proof}
Let $B$ be the vertex matrix of $E$ and $C$ the vertex matrix of $F$.
Without loss of generality, assume $B_1^T = E_{k_1} + \cdots + E_{k_r} +
B_{m_1}^T + \cdots + B_{m_s}^T$.  

We again define a map $\phi: E_1 \rightarrow F^1 \cup F^2$ by 
$$
\phi(e^{ij}) = \left\{ \begin{array}{ll}
       f^{im}f^{m1} & \hbox{if $j = 1$ and $i \not\in K$}
\\
       f^{ij}       & \hbox{else}
                       \end{array}
               \right.
$$
\noindent and extend it to a map $\phi: \mathcal{G}_E \rightarrow
\mathcal{G}_F$.  By arguments similar to those in the proof of
Proposition~\ref{prop-main}, $\phi$ is an injective groupoid homomorphism.
In this case, however, it fails to be onto, and this is because of the
difference between Lemma~\ref{lem-newoldhasinverse} and
Lemma~\ref{lem-oldnewhasinverse}.  Thus it is necessary to point.

It is not hard to see that $\phi$ restricts to an injective groupoid
homomorphism from $\mathcal{G}_{(E,\{v\})}$ to
$\mathcal{G}_{(F,\{v\})}$ which, of course, we also denote by
$\phi$. 

We claim that any finite or infinite path whose first edge
is not new can be
pulled back through $\phi$.  First, find all the new edges in the path.
These are preceded by edges
which are not new, and these (not new)-new pairs can be pulled back.  The
remaining edges are all not new and can be pulled back individually.

Now, fix any $[\alpha, x, \beta]$, with $s(\alpha) = s(\beta) = 1$.
Since $\alpha$ and $\beta$ start at $1$, their first edge cannot be new
(because $1 \not\in M$).  Hence they can be pulled back through $\phi$.
Now, if $x_1$ is not new, $x$ can be pulled back as well, so $[\alpha, x,
\beta]$ has an inverse image.  Since $\phi$ preserves source and range,
this inverse image will be in $\mathcal{G}_{(E,\{v\})}$.  If, on
the
other hand, $x_1$ is new, then we know $x_2$ is not new, so we pull back
the triple $[\alpha x_1, x_2 x_3 \dots , \beta x_1]$.  This shows that
$\phi: \mathcal{G}_{(E,\{v\})} \rightarrow
\mathcal{G}_{(F,\{v\})}$ is onto.  

It is not hard to check that $\phi$ is continuous.  We show that
it is open.  Clearly, $\phi(Z(\alpha, \beta)) \subseteq
Z(\phi(\alpha), \phi(\beta))$.  We show the reverse inclusion:  let
$[\phi(\alpha), y, \phi(\beta)] \in Z(\phi(\alpha), \phi(\beta))$.  If
$y_1$ is not
new, then $y$ has an inverse image and hence $[\phi(\alpha), y,
\phi(\beta)] \in \phi(Z(\alpha, \beta))$.  If, on the other hand, $y_1$ is
a new edge, then note that $y_2$ is not new, so the path 
$y_2y_3\dots$ has an inverse image.  Since 
$$
Z(\phi(\alpha), \phi(\beta)) = \bigcup_{s(f) = r (\alpha)}
Z(\phi(\alpha)f, \phi(\beta)f),
$$ 
$[\phi(\alpha), y, \phi(\beta)] \in \phi(Z(\alpha, \beta))$.
Thus
$\phi(Z(\alpha, \beta)) = Z(\phi(\alpha), \phi(\beta))$, which shows that
$\phi$ is open. 
\end{proof}

\begin{cor}
If $E$ is a graph with no sinks and $F$ is reverse primitively
equivalent to $E$, then $C^*(E)$ is Morita equivalent to $C^*(F)$.  
\end{cor}

\begin{proof} 
If $F$ is a reverse primitive transfer of $E$, then the
result follows easily from the previous proposition and the fact that
pointing a
graph at a cofinal vertex does not change the Morita equivalence class of
its $C^*$-algebra \cite{kprr}.  Since reverse primitive equivalence is
generated by the reverse primitive transfer and permutation, the result
follows.  
\end{proof}

%==========================================================================

\section{Explosions}
\label{sec-explosions}
%==========================================================================

Given a graph $E$, its {\it edge matrix} $A_E$ is an $E^1
\times
E^1$ matrix defined by 
$$
A_E(e,f) = \left\{\begin{array}{ll}
                   1 & \hbox{if $r(e) = s(f)$}\\
                   0 & \hbox{otherwise.}
                 \end{array}  
           \right.
$$ 
The {\it adjoint graph} of $E$ is the graph whose vertex matrix is the
edge matrix of $E$.  In \cite{wat-prim} the explosion of a graph was
defined as a
generalization of the adjoint graph, and it was shown that exploding a
graph does not change its $C^*$-algebra. Since we work with the vertex
matrix instead of the adjacency matrix, our edges are backwards.  To be
consistent with our earlier terminology, we shall call {\it reverse
explosion} what Enomoto, Fujii, and Watatani called explosion, and we
develop a very similar notion, which we shall call {\it explosion}.

Let $E$ be a graph.  Let $v \in E^0$ satisfy $|s^{-1}(v)| > 1$, and fix an 
edge $e$ whose source is $v$.  First assume that $e$
is not a loop, that is,
$v\not=r(e)$. Denote the set of non-loop edges pointing to $v$ by
$K=\{k_1,k_2,\dots\}$ and the set of non-loop edges different from $e$
starting
at $v$ by $M=\{m_1,m_2,\dots\}$. The \emph{edge explosion} $F$ of $E$ at
the edge
$e$ is defined as follows. Split the vertex $v$ into two vertices $v'$ and
$v''$. The source of $e$ is replaced by $v'$. The source of each edge in
$M$ is replaced by $v''$. Every edge in $K$ is replaced by a pair of edges
$k'$, $k''$ having the same source as $k$ and pointing to $v'$ and $v''$
respectively. If there is a loop edge $f$ at $v$ then it is replaced by a loop
$f''$ at $v''$ and an edge $f'$ pointing from $v''$ to $v'$.
The
following picture shows an example of $E$ and its explosion at $e$.
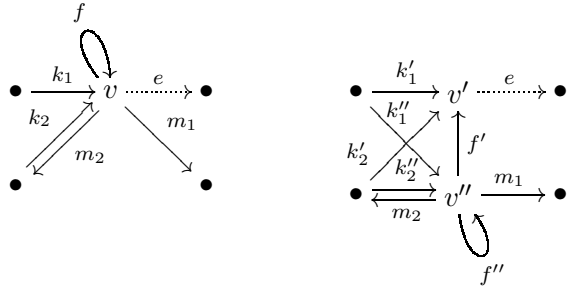
$$ 
\xymatrix{  
\bullet \ar[r]^{k_1}	& v \ar@{.>}[r]^e \ar@<2pt>[dl]^{m_2} \ar[dr]^{m_1} \ar@(lu,u)[]^f & \bullet\\ 
\bullet \ar@<2pt>[ur]^{k_2}	& 	& \bullet 
} 
\qquad\qquad 
\xymatrix{ 
\bullet \ar[r]^{k_1'}_<>(.4){k_1''} \ar[rd]	& {v'} \ar@{.>}[r]^e  & \bullet\\ 
\bullet \ar[ur]^<>(.1){k_2'}	\ar@<2pt>[r]^{k_2''}	& {v''}
\ar@(d,rd)[]_{f''}\ar[u]_{f'} \ar[r]^{m_1}\ar@<2pt>[l]^{m_2}	& \bullet 
} 
$$ 

Next assume that $e$ is a loop. To get the explosion at $e$, split the
vertex
$v$ into $v'$ and $v''$ and change the edges in $M$ and $K$ as before.
Also, replace $e$ by a loop $e'$ at $v'$ and an edge $e''$ pointing from  
$v'$ to $v''$. The following picture shows an example of $E$ and its
explosion at a loop $e$.
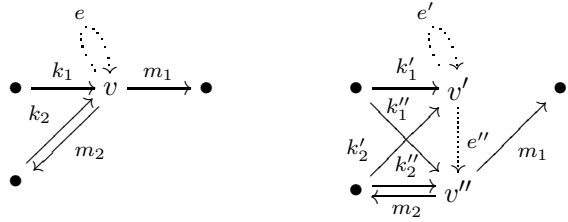
$$ 
\xymatrix{  
\bullet \ar[r]^{k_1}	& v \ar@<2pt>[dl]^{m_2} \ar[r]^{m_1} \ar@{.>}@(lu,u)[]^e & \bullet\\ 
\bullet \ar@<2pt>[ur]	^{k_2}& 	& 
} 
\qquad\qquad 
\xymatrix{  
\bullet \ar[r]^{k_1'}_<>(.4){k_1''} \ar[rd] & {v'} \ar@{.>}[d]^{e''}
\ar@{.>}@(lu,u)[]^{e'} & \bullet\\ 
\bullet \ar[ur]^<>(.1){k_2'}	\ar@<2pt>[r]^{k_2''} & {v''} \ar@<2pt>[l]^{m_2} \ar[ru]_{m_1}	
} 
$$ 

\begin{defn} 
Two graphs $G$ and $E$ are said to
be {\it explosion equivalent} if there is a finite sequence
$E=F_0,F_1,\dots,F_n=G$ of graphs such that for every $i$, either
$F_i$ is an edge explosion of $F_{i+1}$ or $F_{i+1}$ is an edge 
explosion of $F_i$.
\end{defn}  

Consider the following more general notion of explosion, which we call
{\it vertex explosion}.  Fix $v
\in E^0$ with $|s^{-1}(v)| > 1$.  Instead of exploding at an edge whose
source is $v$, we explode at a subset of edges whose source is $v$.  Write
$s^{-1}(v) = M_1 \cup M_2$, where $M_1$ and $M_2$ are disjoint and
nonempty.  Again we split the vertex $v$ into two vertices $v'$ and $v''$
and put an edge from every vertex in $s(r^{-1}(v))$ to both $v'$ and
$v''$.  Also, put an edge from $v'$ to every vertex in $r(M_1)$ and from
$v''$ to every vertex in $r(M_2)$.  If there is a loop edge in $M_1$, we
add an edge from $v'$ to $v''$.  If there is a loop edge in $M_2$, we add
an edge from $v''$ to $v'$.  
If $F$ is an explosion of $E$ at $v$, we always identify $E^0 \setminus \{v\}$
with $F^0 \setminus \{v',v''\}$.  

The following lemma gives a characterization of vertex explosion in terms
of the vertex matrix.  The proof follows immediately from the definition. 

\begin{lem}
\label{lem-explosionmatrix}
Let $F$
be an explosion of $E$ at vertex $v$ with respect to the decomposition
$s^{-1}(v) = M_1 
\cup M_2$.  Let $B$ and $C$ respectively
denote
the vertex matrices of $E$ and $F$.  Then we have the following:
\begin{enumerate}
\item[(i)]  $B_{uw} = C_{uw}$ for every $u \in E^0 \setminus
\{v\}$
and
every $w
\in
F^0 \setminus \{v',v''\}${\rm ;}
\item[(ii)] $C_{wv'} = C_{wv''}$ for every $w \in F^0${\rm ;}     
\item[(iii)] $C_{v''w} = 1 \Leftrightarrow e^{vw} \in M_2$ and $C_{v'w} =
1 \Leftrightarrow e^{vw} \in M_1$ for every $w \in F^0
\setminus \{v', v''\}${\rm .}
\end{enumerate}  
Further, if $B$ and $C$ are 0-1 matrices satisfying {\rm (i)-(iii)},
then the graph of $C$ {\rm(}i.e. the graph whose vertex
matrix is $C${\rm)} is an explosion at vertex $v$ of the graph of $B$.  
\end{lem}

\begin{defn}
Let $E$ be a graph and suppose that $v \in E^0$ satisfies $|s^{-1}(v)| = k
> 1$.  Order the vertices of $E$ so that $v$ is vertex 1, and let
$B^{(1)}$ denote the vertex matrix of $E$.  Now, for $m=1,2,\dots,k-1$,
perform the following procedure.  First find the largest $j$ such that
$B^{(m)}_{1j} = 1$.  Next, insert the row $E_j$ between rows $1$ and $2$
of $B^{(m)}$.  Then duplicate the first column of the resulting matrix.
Finally, change the 1 in the $(1,j)$ position to a 0.  Name this new
matrix $B^{(m+1)}$.  The {\it complete explosion of E at v} is defined to
be the graph of the matrix $B^{(k)}$.  
\end{defn}

Note that, by the previous lemma, each of the $k-1$ steps in the above
procedure corresponds to an explosion at an edge.  We now 
offer the following example to guide the reader through the definition
of complete explosion.

\begin{ex} The graph on the left can be completely exploded in two steps
by exploding at the dotted edge each time. 
The dotted edge becomes the dashed edge in the exploded graph at each stage.
In the matrices, $*$ denotes the unaffected parts of the matrix. 
$$ 
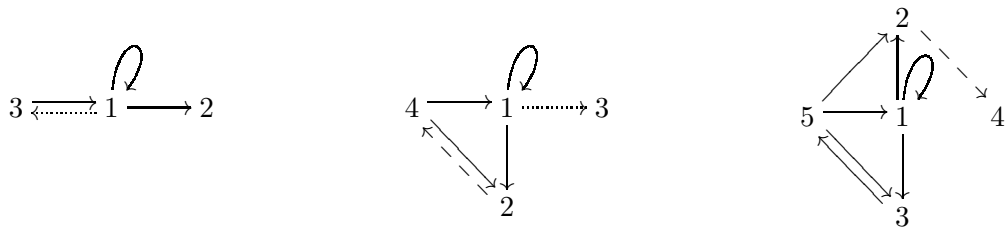
\xymatrix{  
\\
3 \ar@<2pt>[r] & 1 \ar@(u,ru)[] \ar[r] \ar@{.>}@<2pt>[l] & 2
}
\qquad\qquad\qquad 
\xymatrix{
\\  
4 \ar@<2pt>[r] \ar@<2pt>[rd]	& 1 \ar@(u,ru)[] \ar@{.>}[r] \ar[d] & 3 \\
               			& 2 \ar@{-->}@<2pt>[ul]   
} 
\qquad\qquad\qquad
\xymatrix{
                        			& 2 \ar@{-->}@<2pt>[dr] \\
5 \ar@<2pt>[r] \ar@<2pt>[rd] \ar[ur]		& 1 \ar@(u,ru)[] \ar@<2pt>[u] \ar[d] & 4 \\
               					& 3 \ar@<2pt>[ul]   
} 
$$

$$
B^{(1)}=\left(\begin{matrix}          
1 & 1 & 1 \\
0 & * & * \\
1 & * & *
\end{matrix}\right)
\qquad
B^{(2)}=\left(\begin{matrix}          
1 & 1 & 1 & 0 \\
0 & 0 & 0 & 1 \\
0 & 0 & * & * \\
1 & 1 & * & *
\end{matrix}\right)
\qquad
B^{(3)}=\left(\begin{matrix}          
1 & 1 & 1 & 0 & 0 \\
0 & 0 & 0 & 1 & 0 \\
0 & 0 & 0 & 0 & 1 \\
0 & 0 & 0 & * & * \\
1 & 1 & 1 & * & *
\end{matrix}\right)
$$
\end{ex}

The reader may check that if $E$ is a graph with no sinks, the complete
explosion of $E$ at every vertex with more than one edge emanating from it
yields the adjoint graph of $E$.  

\begin{lem}
\label{lem-samerel}
Edge explosion and vertex explosion generate the same equivalence
relation.
\end{lem}

\begin{proof}
Since edge explosion is a special case of vertex explosion, it suffices to
show that an arbitrary graph $E$ and any vertex explosion $F$ of $E$ can
be edge exploded into a common graph.  

Now if $E$ is any graph and $F$ is a vertex explosion of $E$ at $v$, the
reader may verify by examining the vertex matrices that the complete
explosion of $E$ at vertex $v$ coincides with the graph obtained by
performing a complete explosion of $F$ at $v'$ and $v''$.  
\end{proof}

The following is closely related to the notion of explosion defined in
\cite{wat-prim}.

\begin{defn} 
A graph $F$ is a \emph{reverse explosion} of the graph $E$ at a vertex
$v$ if the reverse graph of $F$ is the explosion of the reverse graph of
$E$ at $v$ and $v$ is cofinal.  Two graphs are said to be {\it
reverse explosion equivalent} if there is a finite sequence of reverse
explosions connecting them.
\end{defn}

\begin{prop}  
\label{prop-expleq}
If $E$ is a row-finite graph and $F$ is an explosion of
$E$, then the groupoids of $E$ and $F$ are
isomorphic.   
\end{prop}

\begin{proof} By Lemma~\ref{lem-samerel}, it suffices to prove the case
when $F$ is the explosion of
$E$ at
an edge $e$. First we assume that $e$ is not a loop edge and, we use the
notation $f$, $M$ and $K$ as in the definition of explosion. We define a
map
$\phi:E^\infty\to F^\infty$ as follows. If $x\in E^\infty$ then $\phi$ makes
the following replacements on path segments of $x$:
\begin{align*}
\phi(\cdots ke\cdots)&=\cdots k'e\cdots \\
\phi(\cdots k\overbrace{f\cdots ff}^{n+1} e\cdots)&=
 \cdots k''\overbrace{f''\cdots f''}^n f'e\cdots\\
\phi(\cdots k\overbrace{f\cdots f}^{n} m\cdots)&=
 \cdots k''\overbrace{f''\cdots f''}^n m\cdots \\
\phi(\overbrace{f\cdots ff}^{n+1} e\cdots)&=
 \overbrace{f''\cdots f''}^n f'e\cdots\\
\phi(\overbrace{f\cdots f}^{n} m\cdots)&=
 \overbrace{f''\cdots f''}^n m\cdots\\
\phi(\cdots k\overbrace{fff\cdots }^{\hbox{\tiny{all} $f$}} ) &=
 \cdots k''\overbrace{f''f''f''\cdots}^{\hbox{\tiny{all} $f''$}} \\
\phi(\overbrace{fff \cdots}^{\hbox{\tiny{all} $f$}} )&=
 \overbrace{f''f''f''\cdots}^{\hbox{\tiny{all} $f''$}} 
  \end{align*}
where $k\in K$, $m\in M$ and $n$ is a non-negative integer. It is easy but tedious to check
that $\phi:\mathcal G_E\to \mathcal G_F$ defined by
$\phi(x,k,y)=(\phi(x),k,\phi(y))$ is a bijective homomorphism. 

It remains to check that it is open
and continuous. First we extend $\phi$ to a subset of $E^*$.  If
$\alpha\in E^*$ and $r(\alpha) \neq v$, 
then we define $\phi(\alpha)$ similarly to the above. To show that $\phi$
is
open it
suffices to check that $\phi(Z(\alpha,\beta))$ is open.  First suppose
that $r(\alpha)=r(\beta) \neq v$.  In this case,   
\begin{align*}
\phi(Z(\alpha,\beta)) &= 
\{(\phi(\alpha x,|\alpha|-|\beta|,\beta x) : x\in E^\infty, s(x)=r(\alpha) \}\\
&=
\{(\phi(\alpha)\phi(x),|\phi(\alpha)|-|\phi(\beta)|,\phi(\beta)\phi(x)) : x\in E^\infty, s(x)=r(\alpha) \}\\
&=
\{(\phi(\alpha)y,|\phi(\alpha)|-|\phi(\beta)|,\phi(\beta)y) : y\in F^\infty,
s(y)=r(\phi(\alpha)) \}\\
&=
Z(\phi(\alpha), \phi(\beta)).
\end{align*}  
Now, if $r(\alpha) = r(\beta) = v$, we have three cases.  If $\alpha =
\alpha'k$ and $\beta = \beta'l$ for some $k,l \in K$, then
$\phi(Z(\alpha, \beta)) = Z(\phi(\alpha')k', \phi(\beta')l') \cup Z(\phi(\alpha')k'',
\phi(\beta')l'')$.  If $\alpha = \alpha'k$ and $\beta = \beta'lff\cdots f$, then
$\phi(Z(\alpha, \beta)) = Z(\phi(\alpha')k'', \phi(\beta')l''f''f''\cdots f'')$.  Finally,
in the case where $\alpha = \alpha'kff\cdots f$ and $\beta =
\beta'lff\cdots f$, we have 
$\phi(Z(\alpha, \beta)) = Z(\phi(\alpha')k''f''f''\cdots f'',
\phi(\beta')l''f''f''\cdots f'')$.
Continuity of $\phi$ follows from a similar argument.

The case when $e$ is a loop edge handled similarly, using a slightly different
definition for $\phi$. It now makes the following replacements on $x\in
E^\infty$:
\begin{align*}
\phi(\cdots km\cdots)&=\cdots k''m\cdots \\
\phi(\cdots k\overbrace{e\cdots ee}^{n+1} m\cdots)&=
 \cdots k'\overbrace{e'\cdots e'}^n e''m\cdots \\
 \phi(\overbrace{e\cdots ee}^{n+1} m\cdots)&=
 \overbrace{e'\cdots e'}^n e''m\cdots \\
 \phi(\cdots k \overbrace{ee\cdots}^{\hbox{\tiny{all} $e$}}) &=
 \cdots k \overbrace{e'e'\cdots}^{\hbox{\tiny{all} $e'$}} \\
 \phi(\overbrace{ee\cdots}^{\hbox{\tiny{all} $e$}}) &=
 \overbrace{e'e'\cdots}^{\hbox{\tiny{all} $e'$}}
\end{align*}
where $k\in K$, $m\in M$ and $n$ is a non-negative integer.
\end{proof}

\begin{cor}
\label{cor-expleq}
If $E$ is a graph and $F$ is explosion equivalent to $E$, then
$C^*(E) \cong C^*(F)$.
\end{cor}

\begin{proof}
If $E$ has no sinks, this follows immediately from the proposition.  If
$E$ has sinks, then one can readily verify that $\tilde{F}$ is an
explosion of $\tilde{E}$ and so $C^*(E) \cong C^*(\mathcal{G}_{\tilde{E}})
\cong C^*(\mathcal{G}_{\tilde{F}}) \cong C^*(F)$.
\end{proof}

The following two results are proven similarly to
Proposition~\ref{prop-expleq} and Corollary~\ref{cor-expleq}.

\begin{prop} If $E$ is a row-finite graph and $F$ is the reverse
explosion of $E$ at an edge whose source is $v$, then the groupoid
of $E$
pointed at $v$ and the groupoid of $F$ pointed at $v''$ are isomorphic.
\end{prop}

\begin{cor} If $E$ is a row-finite graph with no sinks and $F$ is a
reverse explosion of $E$ then $C^*(E)$ and $C^*(F)$ are Morita equivalent. 
\end{cor}

%======================================================================

\section{Elementary strong shift equivalence}
\label{sec-SSE}
%======================================================================

A matrix $A$ is \emph{elementary strong shift equivalent} to a matrix $B$
if there
are matrices $R$ and $S$ such that $A=RS$ and $B=SR$ \cite{williams}.
Note that, for any
permutation matrix $P$, $PBP^{-1}$ and $B$ are elementary strong shift
equivalent via
$R = PB$, $S=P^{-1}$.  Thus, any two vertex matrices of the
same graph are elementary strong shift equivalent. Two graphs $E$ and $F$
are said to be elementary strong shift equivalent if
their vertex matrices $B_E$ and $B_F$ are elementary strong shift
equivalent. 

Note that elementary strong shift equivalence is not an equivalence relation. 
The equivalence relation generated by elementary strong shift equivalence is
called {\it strong shift equivalence}. 

Note that if $B_E=RS$ and $B_F=SR$ then the rows and columns of $R$
can be indexed by $E^0$ and $F^0$ respectively. Also the rows and columns of
$S$ can be indexed by $F^0$ and $E^0$ respectively. Using this property we 
define a bipartite \emph{imprimitivity graph} $X$ as follows:  the
set of vertices of $X$ is
the disjoint union of $E^0$ and $F^0$ and the vertex matrix of $X$ is
$$
B_X=\left(\begin{matrix}          
0 & R \\
S & 0
\end{matrix}\right).
$$
The construction of $X$ is due to Ashton \cite{ashton}.

\begin{ex}
Using 
$
R=\left(\begin{smallmatrix}  
1 & 1 & 0 \\
0 & 0 & 1
\end{smallmatrix}\right)
$ and 
$
S=\left(\begin{smallmatrix}  
1 & 0 \\
0 & 1 \\
0 & 1
\end{smallmatrix}\right)
$
we have $E$, $X$ and $F$ 
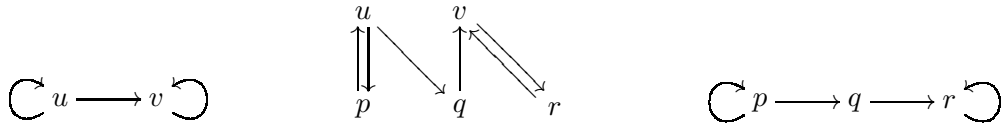
$$ 
\xymatrix{  
\\
u \ar@(dl,ul)[]	\ar[r] & v \ar@(rd,ru)[]
}
\qquad\qquad\qquad 
\xymatrix{  
u \ar@<2pt>[d]\ar[dr]	& v \ar@<2pt>[dr] \\
p \ar@<2pt>[u]			& q \ar[u] & r \ar@<2pt>[ul]
} 
\qquad\qquad\qquad
\xymatrix{  
\\
p \ar@(dl,ul)[]	\ar[r] & q \ar[r] & r \ar@(rd,ru)[]
} 
$$
with vertex matrices
$$
B_E=\left(\begin{matrix}  
1 & 1 \\
0 & 1
\end{matrix}\right)
\qquad\qquad
B_X=\left(\begin{matrix}  
0 & 0 & 1 & 1 & 0 \\
0 & 0 & 0 & 0 & 1 \\
1 & 0 & 0 & 0 & 0 \\
0 & 1 & 0 & 0 & 0 \\
0 & 1 & 0 & 0 & 0 
\end{matrix}\right)
\qquad\qquad
B_F=\left(\begin{matrix}  
1 & 1 & 0 \\
0 & 0 & 1 \\
0 & 0 & 1
\end{matrix}\right)
$$
\end{ex}

\begin{prop}
If $E$ and $F$ are elementary strong shift equivalent, row-finite graphs
and $X$ is the imprimitivity graph, then the
groupoid of
$X$ pointed at $E^0$ is isomorphic to $\mathcal G_E$  and the groupoid of $X$
pointed at $F^0$ is isomorphic to $\mathcal G_F$.   
\end{prop}

\begin{proof}First note that $E^0$ and $F^0$ are automatically cofinal pointing sets.
By symmetry it suffices to show that $\mathcal{G}_{(X,E^0)} \cong \mathcal G_E$. By the
construction of $X$ we have a unique bijection $\phi:E^1\to X^2$ such that
$s=s\circ \phi$ and $r=r\circ \phi$. We extend $\phi$ to $E^*\cup
E^\infty$. It is easy to check that $\phi:\mathcal G_E\to \mathcal
G_{(X,E^0)}$ defined by 
$$
\phi[\alpha,x,\beta]=[\phi(\alpha),\phi(x),\phi(\beta)] 
$$ 
is an isomorphism.
\end{proof}

The $C^*$-algebras of strong shift equivalent graphs are not necessarily
isomorphic, (see Example~\ref{ex-counter1}), but we have:

\begin{cor} If $E$ and $F$ are strong shift equivalent, row-finite
graphs with no sinks then $C^*(E)$ and $C^*(F)$ are Morita equivalent.  
\end{cor}

Following Ashton \cite{ashton}, we say a 0-1 matrix is {\it column
subdivision} if each
of its columns
contains at most one 1.  Two matrices $A$ and $B$ are said to be
elementary strong
shift equivalent with column subdivision if $A = RS$, $B=SR$, and either 
$R$ or $S$ is column subdivision.  Likewise, two graphs are said to be
elementary strong shift equivalent with {\it column subdivision} if their
vertex matrices are.  It was shown in both \cite{ashton} and
\cite{wat-prim} that if two finite 0-1 matrices $A$ and 
$B$ are elementary strong shift equivalent with column subdivision, then 
$\mathcal{O}_A \cong \mathcal{O}_B$.  The following
result, combined with Corollary~\ref{cor-expleq}, provides an alternate
proof of a special case of this fact:

\begin{prop}
Let $E$ and $F$ be graphs with no sinks.  Suppose that $E$ and
$F$ have $n$ and $n+1$ vertices
respectively.  Then $E$ and $F$ are elementary strong shift equivalent
with column
subdivision if
and only if $F$ is an explosion of $E$.  
\end{prop}

\begin{proof} 
Denote the vertex matrix of $E$ by $B$ and the vertex matrix of $F$ by
$C$.  Now suppose that $F$ is an explosion of $E$, with vertex
$v$ splitting into $v'$ and $v''$.  Without loss of generality, we may
assume that $v$ is vertex 1 in $E$ and $v'$ and $v''$ are the first two
vertices of $F$.  That is, the first row and column of $B$ corresponds to
$v$ and the first two rows and columns of $C$ correspond to $v'$ and
$v''$.   
Define $S$ to be $C$ with the first column deleted.  That is, $S_{ij} =
C_{ij}$ for $i=1,\dots,n+1$, $j=1,\dots,n$.  If $R$ is the following $n
\times (n+1)$ column subdivision matrix
$$
R=\left(\begin{array}{ccccc}
1 & 1 & 0 & 0 & 0 \\
0 & 0 & 1 & 0 & 0 \\
 &  &  & \ddots &  \\
0 & 0 & 0 & 0 & 1
\end{array}\right)
$$
\noindent then $B=RS$ and $C=SR$.

Now suppose that $E$ and $F$ are elementary strong shift equivalent with
column
subdivision.  That is, for any choices $B$ and $C$ of vertex matrices of
$E$ and $F$, there exist $R,S$ such that $B = RS$, $C = SR$, and either
$S$ or $R$ is column
subdivision.  Now, $S$ is an $(n+1) \times n$ matrix, so in order for it
to be column subdivision, it must have a zero row.  But a zero row in $S$
would yield a zero row in $C$, and hence a sink in $F$.  Thus it must be
$R$ which is column subdivision.  

Since $R$ is an $n \times (n+1)$ matrix which is column subdivision and
has no zero rows, there exist an $n \times n$ permutation matrix $P$ and
an $(n+1) \times (n+1)$ permutation matrix $Q$ such that 
$$
PRQ=\left(\begin{array}{ccccc}
1 & 1 & 0 & 0 & 0 \\
0 & 0 & 1 & 0 & 0 \\
 &  &  & \ddots &  \\
0 & 0 & 0 & 0 & 1
\end{array}\right)\hbox{.} 
$$
\noindent So by replacing $B$ with $PBP^{-1}$ and $C$ with $QCQ^{-1}$, we
may assume that $R$ has this form.  

Now, if we index the rows of $R$ and the columns of $S$ by
$\{1,2,\dots,n\}$ and the columns of $R$ and the rows of $S$ by
$\{0,1,\dots,n\}$, then we have the following:
\begin{enumerate}
\item[(i)] $B_{ij} = C_{ij}$ for $i=2,\dots,n$, $j=1,\dots,n$; 
\item[(ii)] $C_{i0} = C_{i1}$ for $i=0,\dots,n$;
\item[(iii)]  $C_{0j} + C_{1j} \leq 1$ for $j=0\dots,n$.
\end{enumerate}
\noindent (i) and (ii) are easy to check.  To see (iii), suppose that for
some $j$, $C_{0j}$ and $C_{1j}$ are both 1.  If $j > 1$, then 
$$2=C_{0j}+C_{1j}=\sum_k S_{0k}R_{kj} + \sum_l S_{1l}R_{lj} =
S_{0j}+S_{1j} = \sum_m R_{1m}S_{mj} = B_{1j},$$
which is a contradiction.  If, on the other hand, $j \leq 1$, then similar
calculations show that $B_{11} = 2$.  

These three facts imply that $B$ and $C$ satisfy the three conditions of
Lemma~\ref{lem-explosionmatrix} for some suitable choice of $v$, $M_1$ and
$M_2$.   Thus $E$ is an explosion of $F$.  
\end{proof}

\begin{rem}
Note that the restriction on sinks is necessary in the preceding
proposition, since if 
$$
B=\left(\begin{array}{cc}
1 & 1\\
0 & 0\\
\end{array}\right)\hbox{ and }
C=\left(\begin{array}{ccc}
1 & 1 & 1\\
0 & 0 & 0\\
0 & 0 & 0\\
\end{array}\right)\hbox{,}
$$
\noindent then $B$ and $C$ are elementary strong shift equivalent with
column
subdivision via 
$$
R=\left(\begin{array}{ccc}
1 & 1 & 1\\
0 & 0 & 0\\
\end{array}\right)\hbox{ and }
S=\left(\begin{array}{ccc}
1 & 0\\
0 & 1\\
0 & 0\\
\end{array}\right)\hbox{,}
$$
\noindent but the graph of $C$ is not an explosion of the graph of $B$.
\end{rem}

%------------------------------------------------------------------------
\section{Counterexamples}
\label{sec-counterexamples}
%-------------------------------------------------------------------------
In this section we collect several examples which show that neither
primitive equivalence nor reverse primitive equivalence is implied by
any of the other equivalence relations discussed in this paper.

\begin{ex} 
\label{ex-counter1}
Elementary strong shift equivalence does not imply primitive
equivalence.  This is trivially true because two matrices which are
primitively equivalent must be the same size, while elementary strong
shift
equivalence may change the size.  But this example, taken from
\cite{wat-prim}, shows that elementary strong shift equivalent matrices
need not be
primitively equivalent, even if they have the same size:  if 
$$
R=\left(\begin{array}{ccc}
1 & 0 & 0\\
1 & 0 & 0\\
0 & 1 & 1\\
\end{array}\right)\hbox{ and }
S=\left(\begin{array}{ccc}
0 & 0 & 1\\
1 & 0 & 0\\
0 & 1 & 1\\
\end{array}\right)\hbox{,}
$$ 
\noindent then $\mathcal{O}_{RS} \cong \mathcal{O}_3$ and   
$\mathcal{O}_{SR} \cong \mathcal{O}_3 \otimes M_2$ are not
isomorphic \cite{paschke-salinas}.
Hence the graph corresponding to $RS$ cannot be primitively equivalent to
the graph corresponding to $SR$.  This example also answers negatively a
question
posed in \cite{ashton}, namely, do elementary strong shift equivalent
graphs always yield isomorphic algebras?  
\end{ex}

\begin{ex} Elementary strong shift equivalence does not imply reverse
primitive
equivalence.  Using $R$ and $S$ from the previous example, $S^TR^T =
(RS)^T$ and $R^TS^T = (SR)^T$ are not primitively
equivalent.  This can be verified from the table on page 450 of
\cite{wat-prim}.  We remark in passing that the second graph in the final
row of that table is misprinted.  There should not be a loop on the top
vertex, and there should be a loop added to the lower left vertex. 
\end{ex}

\begin{ex} Reverse primitive equivalence does not imply primitive
equivalence.  Example~\ref{ex-rpe} shows two matrices which are reverse
primitively equivalent, but whose Cuntz-Krieger algebras are
not isomorphic.  Hence they are not primitively equivalent.
\end{ex}

\begin{ex} Explosion equivalence does not imply primitive
equivalence and reverse explosion equivalence does not imply reverse
primitive equivalence.
Consider the matrices
$$
B=\left(\begin{array}{cccc}
1 & 1 & 0 & 1\\
0 & 0 & 1 & 0\\
1 & 1 & 1 & 0\\
1 & 1 & 0 & 1\\
\end{array}\right)\hbox{ and }
C=\left(\begin{array}{cccc}
1 & 1 & 0 & 0\\
0 & 0 & 1 & 1\\
1 & 1 & 1 & 0\\
1 & 1 & 0 & 1\\
\end{array}\right)\hbox{.}
$$ 
\noindent Both are explosions of 
$A=
\left(\begin{smallmatrix}
1 & 1 & 1 \\
1 & 1 & 0 \\
1 & 0 & 1
\end{smallmatrix}\right)
$
so they are explosion equivalent.  But they are not primitively
equivalent.  The reader with a spare afternoon may verify this by checking
that there are 60 elements in the primitive equivalence class of $C$, and
$B$ is not one of them.  Also note that, since $A$ is irreducible, every
vertex is cofinal.  Thus $B^T$ and $C^T$ are reverse explosion
equivalent, but not reverse primitively equivalent.
\end{ex}

\begin{ex} Reverse explosion equivalence does not imply primitive
equivalence and explosion equivalence does not imply reverse primitive
equivalence.  If
$$
B=\left(\begin{array}{ccccc}
0 & 0 & 0 & 1 & 0\\
0 & 1 & 1 & 0 & 0\\
1 & 0 & 0 & 0 & 1\\
0 & 1 & 0 & 0 & 0\\
0 & 1 & 0 & 0 & 0\\
\end{array}\right)\hbox{ and }
C=\left(\begin{array}{ccccc}
0 & 0 & 0 & 0 & 1\\
0 & 1 & 0 & 1 & 0\\
0 & 1 & 0 & 1 & 0\\
1 & 0 & 0 & 0 & 1\\
0 & 0 & 1 & 0 & 0\\ 
\end{array}\right)\hbox{,}
$$ 
\noindent $B$ and $C$ are reverse explosion equivalent because their
transposes are explosions of 
$A=
\left(\begin{smallmatrix}
0 & 0 & 0 & 1\\
0 & 1 & 0 & 1\\
0 & 1 & 0 & 0\\
1 & 0 & 1 & 0
\end{smallmatrix}\right)
$
but they are not primitively equivalent.  Unfortunately, it
requires a computer program to verify this (the primitive equivalence
class of $C$ has 183204 elements), and we could not find a more
manageable example.  
\end{ex}

\begin{ex} Primitive equivalence does not imply reverse primitive
equivalence.  It can be verified using Section 4 of \cite{wat-prim} that 
$$
B=\left(\begin{array}{ccc}
1 & 1 & 1\\
1 & 0 & 0\\
1 & 0 & 0\\
\end{array}\right)\hbox{ and }
C=\left(\begin{array}{ccc}
1 & 1 & 0\\
1 & 0 & 1\\
0 & 1 & 0\\
\end{array}\right)\hbox{}
$$ 
\noindent are primitively equivalent, but not reverse primitively
equivalent.
\end{ex}

It is an open question whether or not primitive equivalence and explosion
equivalence together are enough to characterize graph groupoid
isomorphism.  That is, given two graphs with isomorphic groupoids, are
they explosion-primitive equivalent?  There
are difficulties on both ends.  In particular, given two graphs, determining
whether or not their groupoids are isomorphic is highly non-trivial.  Also,
we currently do not have an efficient algorithm for determining whether
or not two graphs are explosion-primitive equivalent.  There are two
difficulties here.  First, even in the $5 \times 5$ case, some matrices
have primitive equivalence classes with an unmanageable number of elements, so
the computer time required  to check whether two matrices are primitively
equivalent becomes an issue.  Second, we cannot say with
certainty that two matrices are not explosion equivalent.  For example,
consider the matrices $A$ and $B$ from Section~\ref{sec-classification}.  We
have checked that no explosion of $A$ to a $5 \times 5$ matrix is primitively
equivalent to any explosion of $B$ to a $5 \times 5$ matrix.  However,
it may be possible, for example, that $A$ and $B$ can be exploded into $6
\times 6$ (or larger) matrices which are primitively equivalent. 

It would be desirable to have more graph transformations that preserve the
isomorphism class of the groupoid or the $C^*$-algebra, or the Morita
equivalence class of the $C^*$-algebra, of the graph.  Having more
operations would increase our chances of finding a canonical form.

%----------------------------------------------------------------------------
%References
%---------------------------------------------------------------------------

\bibliographystyle{amsplain}

\begin{thebibliography}{1}

\bibitem{ashton}
B.~Ashton, \emph{Morita equivalence of graph ${C}^*$-algebras}, Honors thesis,
  Univ. of Newcastle (1996).

\bibitem{CK1}
J.~Cuntz and W.~Krieger, \emph{A class of ${C}^*$-algebras and topological
  {M}arkov chains}, Invent. Math. \textbf{56} (1980), 251--268.

\bibitem{wat-prim}
M.~Enomoto, M.~Fujii, and Y.~Watatani, \emph{${K}_0$-groups and classifications
  of {C}untz-{K}rieger algebras}, Math. Japon. \textbf{26} (1981), 443--460.

\bibitem{franks}
J.~Franks, \emph{Flow equivalence of subshifts of finite type}, Ergo. Th. \&
  Dynam. Sys. \textbf{4} (1984), 53--66.

\bibitem{kpr}
A.~Kumjian, D.~Pask, and I.~Raeburn, \emph{{C}untz-{K}rieger algebras of
  directed graphs}, Pacific J. Math (to appear).

\bibitem{kprr}
A.~Kumjian, D.~Pask, I.~Raeburn, and J.~Renault, \emph{Graphs, groupoids, and
  {C}untz-{K}rieger algebras}, J. Funct. Anal. \textbf{144} (1997), 505--541.

\bibitem{paschke-salinas}
W.L.~Paschke and N.~Salinas, \emph{Matrix algebras over $\mathcal{O}_n$},
  Mich. Math. J. \textbf{26} (1979), 3--12.

\bibitem{pask/rae-K}
D.~Pask and I.~Raeburn, \emph{On the {K}-theory of {C}untz-{K}rieger algebras},
  Publ. RIMS Kyoto Univ. \textbf{32} (1996), 415--443.

\bibitem{rordam-K}
M.~R\o rdam, \emph{Classification of {C}untz-{K}rieger algebras},
  {K}-theory \textbf{9} (1995), 31--58.

\bibitem{williams}
R.~Williams, \emph{Classification of subshifts of finite type}, Ann. of
  Math. \textbf{98} (1973), 120--153.  Errata \textbf{99} (1974),
  380--381.

\end{thebibliography}

\providecommand{\bysame}{\leavevmode\hbox to3em{\hrulefill}\thinspace}

\vspace*{.25in}
\end{document}